\documentclass[journal]{IEEEtran}

\ifCLASSINFOpdf
\else
\fi

\usepackage[cmex10]{amsmath}
\usepackage{amsmath, amssymb}

\usepackage{ dsfont }
\usepackage{url}
\usepackage{cancel}
\usepackage{ mathrsfs }
\usepackage{theorem}
\usepackage{graphicx}

\usepackage{soul,xcolor}

\usepackage{todonotes}

\DeclareFontEncoding{LGR}{}{}
\DeclareTextSymbol{\~}{LGR}{126}

\theoremstyle{remark}
\theoremheaderfont{\normalfont\bfseries}
\theorembodyfont{\normalfont}
\newtheorem{theorem}{Theorem}[section]

\newtheorem{definition}[theorem]{Definition}

\newtheorem{remark}[theorem]{Remark}

\usepackage[short]{optidef}

\usepackage{stackengine}
\usepackage{scalerel}
\usepackage{enumitem}
\usepackage{ gensymb }
\usepackage{subcaption}
\usepackage[numbers,sort&compress]{natbib}
\usepackage{mathtools}
\usepackage{booktabs}
\usepackage{float}

\usepackage[subtle]{savetrees}
\allowdisplaybreaks[4]
\renewcommand{\baselinestretch}{1}
\usepackage[font=small,skip=1.9pt]{caption}
\usepackage[T1]{fontenc}
\usepackage[utf8]{inputenc}
\usepackage[english]{babel}

\renewcommand{\IEEEbibitemsep}{0pt plus 2pt}
\makeatletter
\IEEEtriggercmd{\reset@font\normalfont\fontsize{3.9pt}{3.40pt}\selectfont}
\makeatother
\IEEEtriggeratref{1}
\begin{document}
\bstctlcite{IEEEexample:BSTcontrol}
\setstcolor{red}
%
\title{Optimal corrective dispatch of \\uncertain virtual energy storage systems}

\author{Mahraz~Amini,~\IEEEmembership{Student Member,~IEEE,}
       and Mads~Almassalkhi,~\IEEEmembership{Senior Member,~IEEE,}
\thanks{The authors are with the Department of Electrical and Biomedical Engineering at The University of Vermont, Burlington, VT 05405, USA. This work was supported by the U.S. Department of Energy’s Advanced Research Projects Agency - Energy (ARPA-E) award DE-AR0000694.}}



\maketitle

\begin{abstract}
High penetrations of intermittent renewable energy resources in the power system require large balancing reserves for reliable operations. Aggregated and coordinated behind-the-meter loads can provide these fast reserves, but represent energy-constrained and uncertain reserves (in their energy state and capacity). 
To optimally dispatch uncertain, energy-constrained reserves, optimization-based techniques allow one to develop an appropriate trade-off between closed-loop performance and robustness of the dispatch. Therefore, this paper investigates the uncertainty associated with energy-constrained aggregations of flexible, behind-the-meter distributed energy resources (DERs). The uncertainty studied herein is associated with estimating the state of charge and the capacity of an aggregation of DERs (i.e., a virtual energy storage system or VESS). To that effect, a risk-based chance constrained control strategy is developed that optimizes the operational risk of unexpectedly saturating the VESS against deviating generators from their scheduled set-points. The controller coordinates energy-constrained VESSs to minimize unscheduled participation of and overcome ramp-rate limited generators for balancing variability from renewable generation, while taking into account grid conditions. To illustrate the effectiveness of the proposed method, simulation-based analysis is carried out on an augmented IEEE RTS-96 network with uncertain energy resources and temperature-based dynamic line ratings. 

\end{abstract}


\section{Introduction}

Conventional generators, such as fast-ramping gas generators, have provided reliable balancing reserves to meet the variability of traditional demand. However, with the increasing penetration of wind and solar PV generation in power systems, more fast reserves are needed. If these reserves were provided by conventional thermal generation, the generators would be operated at reduced power and can lead to an increase in idling, which is economically inefficient and increases harmful emissions. Rather than managing the net-load variability (i.e., demand minus renewable generation) with conventional generators, a grid-scale energy storage system (ESS) can provide fast-acting reserves from rapid charging or discharging events. However, even though the cost of an ESS has declined significantly in the last decade, ESSs remain expensive options. In addition, less costly alternatives, such as energy services provided by aggregated DERs, are advancing rapidly~\cite{meyn2018compare}.

In fact, flexible, behind-the-meter loads, such as thermostatically controllable loads (TCLs), electric vehicles (EVs), and residential batteries, can be coordinated and aggregated 
to form a \textit{virtual energy storage system} (VESS). A VESS can provide grid services similar those of an ESS, including synthetic (i.e., demand-side) reserves~\cite{Mathieu2015arbitrage, meyn2015ancillary, almassalkhi2017packetized}.
 This is due to the fact that control actions that increase/decrease the power generated by conventional generators can be equivalently provided by a decrease/increase in the aggregate demand. 
While the core concepts underpinning autonomous demand response (ADR) can be traced back to the early 1980s~\cite{morgan1979electric, schweppe1980homeostatic}, the VESS technology available today is still in the early stages, but advancing  rapidly. Recently, researchers have developed generalized energy-based models for aggregating and coordinating DERs that are very similar to that of a classic charge/discharge battery model~\cite{Mathieu2013Arbitrage, Hao:VB2015, Martinez2015balance, Mathieu2015arbitrage, vrakopoulou2017chance, Chakraborty2018VB, Madjidian2018aggregate, Cammardella2018CA}.  The VESS model presented herein takes into account the general, first-order energy dynamics similar to those found in the literature and consider the box constraints on the energy state and (charge/discharge) power dispatch bounds. 

However, a flexible load has its own baseline consumption that is a function of many exogenous and uncertain parameters (e.g. hot water usage, arrival or departure time of EVs, etc.). In addition, individual loads may not be directly controlled by the operator. Thus, unlike a physical grid-scale battery, a VESS's energy state, energy capacity, and other parameters (e.g., power limits) that define its available flexibility are inherently time-varying and uncertain. For example, in ~\cite{zhang:aggDRuncertain2016}, the uncertainty (due to DER failures and repairs) is quantified for a collection of DERs and a flexibility capacity-duration-probability curves is constructed to allow a VESS operator to more effective bid in uncertain capacity. Thus, to benefit the most from the availability of uncertain and energy-constrained VESS-based reserves, careful design of predictive control techniques are necessary to optimize the VESS dispatch.

{\color{black}Since power systems are suffused with constraints and limits, model predictive control (MPC), with the ability to consider multiple inputs and outputs and the temporal coupling inherent to energy storage and dynamic line rating (DLR), makes it a useful tool for corrective dispatch of energy resources. Therefore, MPC has been widely used for set-points optimization of grid resources, including demand and energy storage applications, e.g., please see~\cite{almassalkhi2015modela,almassalkhi2015model,amini2018tradingoff}.} 
 For a general overview of MPC, please see~\cite{mayne2000constrained}. MPC operates the system over a receding horizon by considering forecasts and power and energy states of the available resources. This makes MPC a particularly useful method within which, one can analyze the trade-off between dispatching ramp-rate limited conventional generators and energy-constrained, uncertain VESS resources. Authors in~\cite{Guo2018P1,Guo2018P2} solve a stochastic multi-stage OPF problem that determines power scheduling policies in order to balance operational costs and network constraint violations. However, only the solar PV power injection is assumed to be uncertain, while the dispatchable energy storage devices are deterministic. Unlike that work, we consider dispatchable energy storage resources that are uncertain and develop a tractable approach to solve a stochastic multi-stage OPF problem.
{\color{black}
Indeed, a predictive technique is presented in~\cite{vrakopoulou2017chance,li2017chance}, which develops a reserve scheduling framework that manages uncertain renewable generation and demand-side reserves. They take into account the uncertainty in the capacity of the controllable load over multiple periods.  However, these approaches do not consider energy state uncertainty of the aggregated flexible loads and focus on solving a reserve market clearing problem every 15~minutes with respect to the static power limits of the transmission lines. Herein, we propose to optimize the use of these uncertain VESSs on a minutely timescale as reserves are activated and we seek to deliver the required flexibility to absorb the high level of variability inherent to a sustainable energy future.


Under high penetrations of renewable generation, the resulting power flows make it challenging to 
reliably operate a power system under conventional static thermal (MVA) limits~\cite{sugihara2017evaluation}.  Instead, we consider the dual role of temperature-based dynamic line rating (DLR) on transmission operations as it relates to integrating VESSs (by providing thermal inertia) and increased renewable generation (by temporarily increasing line flow limits). The temperature-based rating is a function of line sagging requirements, which determines the maximum thermal expansion of a conductor (i.e., conductor temperature) to satisfy a minimum ground clearance height~\cite{IEEE738}.  
With DLR, significant operational flexibility can be gained by coupling the line electrical and thermal behaviour and inertia at the cost of additional sensing~\cite{Banakar2005I,Banakar2005II}.  Thus, temperature-based ratings provide the real-time apparent power capacity of lines based on actual and estimated operating conditions, which allows a predictive controller to naturally leverage feedback to ensure reliable operating conditions without having to employ overly conservative limits on the lines' power transfer capabilities. That is,  the relationship between nonlinear AC line flows and uncertain VESSs is challenging to study in the static sense, but with feedback, the DLR formulation can dynamically adapt to these complex interactions, which is exactly how the proposed risk-based VESS dispatch approach complements the overall DLR scheme.

The contributions of this manuscript are the following: $\mathbf{(i})$~as far as the authors are aware,  prior work on stochastic optimal power flow (OPF) methods focuses mainly on the uncertainty of power injections (e.g., wind and demand), which temporally decouples the OPF problem and avoids the challenges of multi-period optimization under uncertainty, e.g.,~\cite{li2016impact,Zhang2017}. Unlike those works, this paper presents a general VESS model with uncertainty in both the estimate of the state of charge and the prediction of energy capacities and incorporates this uncertain VESS into a stochastic, multi-period OPF framework; $\mathbf{(ii})$~by integrating uncertain energy resources with the electro-thermal coordination of line temperature dynamics, we are able to leverage the inertia and the complementary time-scales for control of both energy storage and line limits to effectively manage uncertainty and grid constraints over multiple time-steps; 
and $\mathbf{(iii)}$~with an analytical reformulation, we present a risk-based, chance-constrained MPC (RB-CC-MPC) approach that co-optimizes the delivery of responsive VESS resources against the operational risk inherent to the VESSs' uncertain energy capacities and states of charge. We then show that this risk-based approach consistently outperforms robust and deterministic approaches under different levels of available information on the nature of the VESSs' uncertainty.}

The rest of this paper is organized as follow. In Section~II, we summarize the proposed control framework and discuss the role and interactions of the OPF problem within a reference-tracking, predictive controller. Section~III details the system models while Section~IV describes the nature and management of uncertainty in a chance-constrained formulation. Section~V introduces the novel risk-based chance-constrained approach to manage uncertainty. Via a simulation-based case-study on the IEEE RTS~96 test system augmented with VSSEs, Section~VI illustrates the analytical formulation and compares the proposed risk-based method against deterministic and robust approaches. Section~VII summarizes the key results and describes future work directions.

\section{System Operation and Control}
\begin{figure}[t!]
      \centering
       \includegraphics[width=\columnwidth]{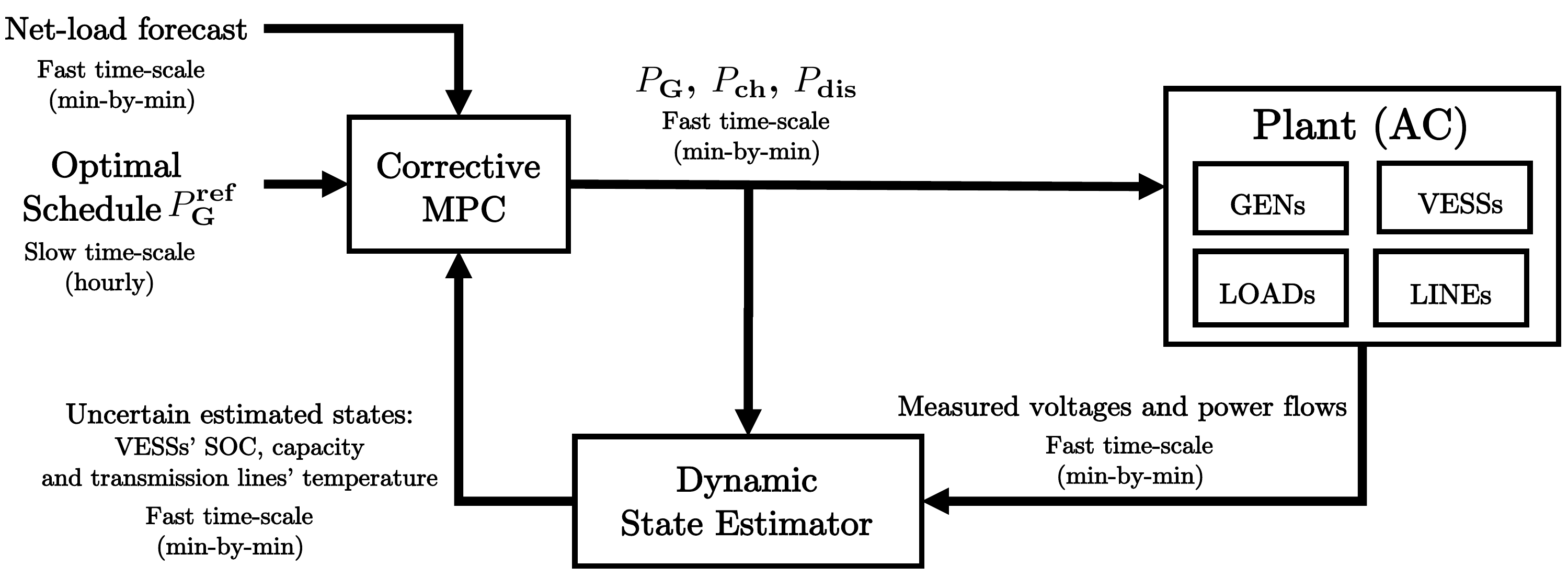}
      \caption{{\small Overview of control scheme showing controller including economic dispatch (slow) and corrective MPC (fast) part and how each part is related to the power grid. }}
      \label{fig:overview}
   \end{figure}

{\color{black} 
Based on updated forecasts of demand and renewable generation, economic dispatch computes a secure and economically optimal schedule for the available generators. However, the uncertainty inherent to solar PV and wind forecasts, as well as uncertainty in electrical demand, results in power imbalances that make previously computed set-points sub-optimal. Rescheduling the generators frequently and significantly accumulates cycling costs and economic penalties to the system operator~\cite{caiso}. With responsive VESS resources and temperature-based DLR, corrective power system operations that leverage feedback represents a valuable and inherently robust and dynamic alternative to conventional spinning reserve. Corrective control refers to the coordination of responsive grid resources immediately \textit{after} a disturbance occurs to drive the system back from an economically sub-optimal or stressed system state to an economically optimal normal operating state~\cite{almassalkhi2015model}. An overview of the proposed uncertainty-aware control strategy is provided in in Fig.~\ref{fig:overview}.

 While the focus of this paper is on the corrective part of the controller, the details of the standard ED can be found in~\cite{wood2012power} and are beyond the scope of this paper. Thus, the output of the market layer ($P_{\text{G}}^{\text{ref}}$ in Fig.~\ref{fig:overview}) satisfies techno-economic objectives such as cost and security. {\color{black}Since the proposed corrective controller's time-step $T_s$ ($\approx1$~minute) is much shorter than the updated of the market-based reference signals coming from ED ($\approx15-60$~minutes), linear interpolation is employed to form the reference trajectory.}
 
 In the faster control layer, the VESSs represent aggregated DERs and provide flexibility in the form of synthetic balancing reserves. Therefore, based on the dynamic states (i.e., power and energy states of VESSs, power states of generators and thermal states of the transmission lines) and forecasts of the system and available resources, the trajectory-tracking MPC produces a corrective dispatch every minute to respond to forecast errors and other disturbances. The MPC minimizes the deviation of generators and flexible loads from the economic reference trajectory while satisfying physical and operational grid constraints. The goal of this manuscript is then to develop and present a risk-based, multi-period OPF formulation that explicitly considers the uncertainty of VESSs in the optimal dispatch. The predictive VESS and OPF models are presented next. 
 

}

 
\subsection{Nomenclature} \textcolor{black}{
The key mathematical symbols are defined next to provide a common reference. However, all mathematical symbols are also described as they are introduced in the remainder of the manuscript:
\begin{itemize}[label={},leftmargin=*]
    \item $c_{G,i}$: cost of deviating from the reference for generator $i$;
    \item $c_{T,ij}$: cost of temperature overload on line $(i,j)$;
    \item $P_{\text{G},i}$: power output of generator $i$;
    \item $R_{G,i}$: maximum ramp rate of generator $i$;
    \item $\Delta P_{G,i}$: change in power output of generator $i$;
    \item $\underline{P_{Gi}}$ ($\overline{P_{Gi}}$): minimum (maximum) production of generator $i$;
    \item $P_{\text{L},n}$: real power consumed by load $n$;
    \item $p_{ij}$: real power flowing on line $(i,j)$;
    \item $p_{ij}^{\text{loss}}$: real power loss on line $(i,j)$;
    \item $\Delta T_{ij}$: temperature overload on line $(i,j)$ with respect to $T_{ij}^{\text{lim}}$;
    \item $P_{\text{ch},i}$ ($P_{\text{dis},i}$): charging (discharging) power of VESS $i$;
    \item $\overline{P_{\text{ch},i}}$ ($\overline{P_{\text{dis},i}}$): maximum charging (discharging) power of VESS $i$;
    \item $\eta_{\text{ch}}$  ($\eta_{\text{dis}}$): charging (discharging) efficiency of VESS $i$;
    \item $S_i$: state of the charge of the VESS $i$;
    \item $\overline{S_i}$ ($\underline{S_i}$): maximum (minimum) energy limit of the VESS $i$;
    \item $R_{\text{ch}} (R_{\text{dis}})$: maximum charging (discharging) ramp rate of VESS $i$.
\end{itemize}
}
 
\color{black}
\section{Predictive model for corrective control}
The corrective MPC scheme is summarized by the following:
\begin{enumerate}
\item At time $k$, with estimates of initial state of the charge (SOC), line temperatures, generator operating points, updated net-load forecasts, and updated generator economic dispatch schedule, the MPC solves a finite-horizon open-loop optimal control problem, over interval [$k,\, k+M$]. This produces a schedule of control actions that describe charging (discharging) rates for { \color{black} VESS} and re-dispatch signals for generators.
\item Apply only the control actions corresponding to time $k$
\item Measure/estimate the system's dynamic states based on the realized demand and renewable generation at time $k+1$.
\item Go to $1)$
\end{enumerate}

The open-loop MPC optimization problem is defined below for a power system network $\mathcal{E}=(\mathcal{N}, \mathcal{L})$ with bus $i \in \mathcal{N}$ and line $(i,j) \in \mathcal{L}$. The sets $\mathcal{G}$ and $\mathcal{V}$ represent conventional generators and VESSs, respectively. The objective function seeks to minimize the deviation of generator outputs from the scheduled set-points while penalizing line temperature overloads as follows:

\small
\begin{mini!}[3]
  {P_\text{G},P_{\text{ch}},P_{\text{dis}}}{ \sum_{l = k}^{k+M} \Bigg( \sum_{\forall i \in N_\text{g}} c_{\text{G},i} \left( P_{\text{G},i}[l]-P_{\text{G},i}^{\text{ref}}[l] \right)^2  + \sum_{\forall ij \in \mathcal{E}} c_{\text{T},ij} {\Delta\hat{\text{T}}_{ij}}^2 \Bigg)}
   {\label{eq:optVESS}}{}
   \addConstraint { \textbf{Power balance: } \forall i \in \mathcal{N} } \nonumber 
   \addConstraint { \sum_{n \in \Omega_i^G}P_{\text{G},n}[l]} {= \sum_{n \in \Omega_i^L} {P_{\text{L},n}[l] + \sum_{n \in \Omega_i^B} P_{\text{ch},n}[l]-P_{\text{dis},n}[l] + \sum_{j \in \Omega_i^N}p_{ij}[l] }
   \label{eq:const2}} 
      \addConstraint { \textbf{Conventional generators: } \forall i \in \mathcal{G}} \nonumber 
    \addConstraint{  P_{\text{G}_i}[l+1]}{ = P_{\text{G}_i}[l] + \Delta P_{\text{G}_i}[l] \label{eq:const5}}
   \addConstraint{  \underline{P_{\text{G}_i}}\leq P_{\text{G},i}[l]}{ \leq \overline{P_{\text{G},i}} \label{eq:const5plus}}
       \addConstraint{  -T_sR_{\text{G},i} \leq \Delta P_{\text{G}_i}[l]}{ \leq T_sR_{\text{G},i} \label{eq:const11}}
          \addConstraint { \textbf{Temperature-based  line rating: } \forall ij \in \mathcal{L}} \nonumber 
      \addConstraint{  p_{ij}[l]}{ = b_{ij}\left(\theta_i[l]-\theta_j[l]\right) + \frac{1}{2} p_{ij,k}^{\text{loss}} \label{eq:line_lineflow} }{}
         \addConstraint{  p_{ij}^{\text{loss}}[l]}{ = R_{ij}\left(b_{ij}\left(\theta_i[l]-\theta_j[l]\right) \right)^2 \label{eq:lineloss_lineflow} }{}
        \addConstraint{  \Delta p_{ij}^{\text{loss}}[l]}{ = p_{ij}^{\text{loss}}[l] - p_{ij,*}^{\text{loss}}  \label{eq:const3_1} }
       \addConstraint{  \Delta T_{ij}[l+1]}{ = \tau_{ij}\Delta T_{ij}[l] + \rho_{ij}\Delta p_{ij}^{\text{loss}}[l] \label{eq:Tline_dyn}}
              \addConstraint{  \Delta \hat{T}_{ij}[l]}{ = \text{max} (0,\Delta T_{ij}[l]) \label{eq:deltaT}}{\mkern-227mu}
   \addConstraint { \textbf{VESSs: } \quad \forall i \in \mathcal{V}} \nonumber 
            \addConstraint{ 0 \leq P_{\text{ch},i}[l]}{ \leq \overline{P_{\text{ch},i}} \label{eq:const6}} 
     \addConstraint{  0 \leq P_{\text{dis},i}[l]}{ \leq \overline{P_{\text{dis},i}} \label{eq:const7}} 
            \addConstraint{ -T_sR_{\text{ch},i} \leq P_{\text{ch},i}[l+1] - P_{\text{ch},i}[l]}{ \leq T_sR_{\text{ch},i} \label{eq:const14}}
                   \addConstraint{ -T_sR_{\text{dis},i} \leq P_{\text{dis},i}[l+1] - P_{\text{dis},i}[l]}{ \leq T_sR_{\text{dis},i} \label{eq:const15}}
      \addConstraint{\underline{S_i} \leq  S_i[l+1]}{ = S_i[l] + T_s\bigg(\eta_{\text{ch},i} P_{\text{ch},i}[l] - \frac{1}{\eta_{\text{dis},i}} P_{\text{dis},i}[l]\bigg) \leq \overline{S_i} \label{eq:const12}}   
    \addConstraint{  S_i[k]}{ = S_{i,k}^{\text{est}} \label{eq:const19}}
\end{mini!}
\normalsize
where $\Omega_i^G$, $\Omega_i^L$, $\Omega_i^B$ and  $\Omega_i^N$ represent set of generators, demands, energy storage devices (VESSs), and neighboring nodes connected to node~$i$, respectively. 
Constraints~\eqref{eq:const2} to~\eqref{eq:const19} must be satisfied for $\forall l \in [k,\ k+M-1]$ where the four groups of constraints in the MPC formulation are described below. Next, we discuss each group.

\subsection{  Power balance in~\eqref{eq:const2}} Based on Kirchhoff's laws, the net power flow into any node must equal the net power flow out. Generators may inject power, $P_G$ and loads may consume power $P_L$ at each node $i$. If VESSs are available at a node, then positive (negative), $P_{\text{ch}}-P_{\text{dis}}$, corresponds to additional consumption (generation).
\subsection{Conventional generators in~\eqref{eq:const5} to~\eqref{eq:const11}:} Each conventional generator is described by its production state, $P_\text{G}$, which must be within generator limits, as shown in~\eqref{eq:const5} and~\eqref{eq:const5plus}. Furthermore, due to the thermal nature of the generators, the ramp rate of generators are limited to their ramp-rate limit, $R_G$,  as shown in~\eqref{eq:const11}. VESSs are particularly helpful to overcome limitations imposed by the ramp-rate limits.

\subsection{Transmission lines in~\eqref{eq:line_lineflow} to~\eqref{eq:deltaT}}

{\color{black} The temperature-based line ratings provide a mechanism through which the uncertain power injections from renewable generation can be absorbed and directed via MPC's feedback and the optimized VESSs dispatch. 
The heat gain of transmission line $(i,j)$ is a function of ohmic losses (i.e. $I_{ij}^2 r_{ij}$). Thus, it is necessary to include line losses in the power flow model. 
Since all values are per-unit (p.u.), and the voltage magnitudes of all buses are close to $1$ p.u., the magnitudes of the current and power flows on the respective lines are approximately equal (i.e., $|I| \approx |S|$ ). Therefore,  line losses can be effectively approximated in proportion to the square of the power flow~\cite{Banakar2005I, almassalkhi2015model}.
In general, the AC power flow between bus $i$ and $j$, $p_{ij}$, is the solution to a set of nonlinear, algebraic equations. 
To ensure a tractable approach at the timescale of interest, a suitable convex relaxation  of the AC power flow equation has been adopted from~\cite{almassalkhi2015modela}.  Since the MPC executes on a fast timescale relative to the VESS and line temperature time-constants and the linearized model is updated via feedback from estimating line losses, $p_{ij,k}^\text{loss}$, line temperatures, and VESS states, the model is sufficiently accurate for control on AC transmission networks.

IEEE Standard~738~\cite{IEEE738} defines the current-temperature relationship of bare overhead conductors and has been employed herein to calculate the conductor temperature. To allow for a tractable implementation of MPC scheme, temperature dynamics of transmission lines are linearized around the equilibrium point $T^* = T^{\text{lim}}$, where $T^{\text{lim}}$ is computed from steady-state conditions with line current at ampacity (i.e., set $p_{ij}^{\text{loss}^*} = (I^{\text{lim}})^2 r_{ij}$, where $r_{ij}$ is the per-unit resistance of line $ij$). The linearized temperature dynamics of the transmission lines are given by~\eqref{eq:Tline_dyn}. Since line losses are approximated in proportion to the square of the power flow~\eqref{eq:lineloss_lineflow}, its respective constraint is non-convex in ${\theta}_{ij}$. Therefore, a convex relaxation is employed (i.e., $p_{ij}^{\text{loss}} \ge R_{ij}b_{ij}^2 (\theta_i - \theta_j)^2$) that is provably binding at optimality for lines that are overloaded since $\Delta \hat T_{ij}$ is in the objective function as detailed in~\cite{almassalkhi2015model}. This achieves the desired measure of control over the line flows.

The MPC scheme computes control actions that drive line temperatures below limits, and, as long as they are below limits, there is no benefit in further reducing line temperatures. The non-convex constraint in~\eqref{eq:deltaT} achieves this purpose.  However, the  non-convex constraint can be relaxed with an equivalent (convex) linear formulation: 
$$ \Delta \hat T_{ij}[l] = \max (0,\Delta T_{ij}[l]) \Longleftrightarrow 0 \leq \Delta \hat{T}_{ij} \phantom{l} \land \phantom{l}  \Delta T_{ij} \leq \Delta \hat{T}_{ij}.$$
All together, this implies that our ``lossy-DC'' grid model formulation is convex. Next, we describe the convex VESS model.

\subsection{Virtual energy storage system in~\eqref{eq:const6} to~\eqref{eq:const19}}

In this paper, responsive VESSs, which are available throughout the network, have a baseline consumption (i.e., aggregated baseline consumption of individual flexible loads in a VESS), and are allocated as balancing reserves. Note also that the VESS models used in this manuscripts are agnostic to the specifics of the coordination scheme. By shifting the VESS's controllable load in time, the VESS can respond to the mismatches caused by forecast errors. Any decrease (increase) in the consumption of the VESS relative to the baseline can be translated as discharging (charging) the VESS. Each VESS is described by an estimated SOC and the amount of power it provides to (consumes from) the grid. At time $k$, the initial SOC of a VESS is given by a dynamic state estimator ~\eqref{eq:const19} and the SOC of a VESS over the prediction horizon is defined by the discrete integrator dynamics as shown in~\eqref{eq:const12}. Non-negative scalar  $P_{\text{ch}}$ ($P_{\text{dis}}$) represents charging (discharging) power of a VESS and the charging (discharging) efficiency is denoted  by $\eta_{\text{ch}}$ ($\eta_{\text{dis}}$).

Charging (discharging) power and SOC of the VESSs are subject to constraints~\eqref{eq:const6},~\eqref{eq:const7} and~\eqref{eq:const12} where $\overline{P_{\text{ch}}}$ ($\overline{P_{\text{dis}}}$) and $\overline{S}$ ($\underline{S}$), respectively, represent maximum charging (discharging) power and the maximum (minimum) energy capacity of VESS.

Since VESSs represent the aggregate effects of coordinated DERs, they inherit the characteristic timescales of the coordination schemes that underpin them. That is, in general, coordination schemes do not offer instant control over all DERs in a fleet, but are subject to separate internal control, actuation, and communication loops~\cite{Duffaut2018tpwrsP1}. These cyber-physical control considerations generalize themselves as 
ramp-rate limits on the charging ($R_{\text{ch}}$) and discharging ($R_{\text{dis}}$) of VESSs as shown in~\eqref{eq:const14} and~\eqref{eq:const15}. At high levels of renewable penetration, since the VESS's are responsive, they represent a valuable resource to overcome demand-supply imbalances. However, unlike a conventional generator, the VESS's energy-constrained characteristics necessitate careful management of its state of charge.

\begin{remark}[Simultaneous charging/discharging]
For most physical ESSs, simultaneous charging and discharging is not physically realizable, which necessitates complex schemes to guarantee an optimal realizable solution~\cite{almassalkhi2015modela}. However, since a VESS coordinates a large population of diverse DERs, a VESS can engender simultaneous charging and discharging commands across the population, e.g., please see~\cite{DuffautEspinosa:2018PSCC, Duffaut2018tpwrsP2}, where residential batteries can strategically discharge to indirectly supply demand for (concurrently charging) hot water heaters that would otherwise become too cold.  
This VESS capability avoids the need for non-convex, complementarity constraints (i.e., $P_{\text{ch},i}[l]P_{\text{dis},i}[l]=0$). 
\end{remark}

Unlike grid-tied batteries, the amount of flexibility available to the system operator is time-varying and uncertain, which is the focus of this manuscript. That is, the flexibility available to the system operator from a VESS can be translated into upper and lower bounds on the VESS's energy state. These upper and lower bounds are functions of different stochastic quantities, such as human behavior and weather. To capture these considerations, VESSs herein are formulated probabilistically and are modeled with chance constraints.

\section{Uncertainty management}

Individual flexible loads are subject to device-specific background effects (such as hot water usage or EV driving patterns), however, a VESS represents a macro-level object. Thus, these time-varying and stochastic background processes are aggregated and realize themselves as uncertainty in the VESS's estimated energy bounds and state. For example, the upper energy capacity limit of the VESS is uncertain and must be estimated and predicted from a separate data-driven model. That is, the flexibility offered by each device is uncertain and represents an independent random variable (i.e., background usage of each device is independent). Therefore, a VESS's energy capacity can be approximated as an aggregation of random variables. Moreover, in contrast with the grid-scale batteries, the actual SOC of the VESSs can not be measured directly and a dynamic state estimation method could be employed to estimate the SOC of a VESS at each time step (e.g. an Extended Kalman Filter, such as in~\cite{espinosaaggregate} for an example of a VESS with state estimation). State estimation of a VESS's SOC is subject to uncertainty inherent in any state estimation method. In addition, due to the nature of Kalman filters, the noise process represents another stochastic effect. 

Note that the exact probability distribution from which the random variables are drawn that describe the energy capacity prediction and SOC estimation errors is generally not known. However, while the distribution may not be known, a VESS operator can leverage available historical and online data to estimate first and second moments of the energy capacity prediction and SOC estimation errors~\cite{zhang:aggDRuncertain2016}. In addition, it is reasonable to consider that the VESS's underlying distribution may satisfy certain conditions, such as being unimodal or having finite support. Without loss of generality, in this manuscript, we consider the case where the VESS operator does not know the exact distribution, but has information that suggests that the distribution is unimodal\footnote{A unimodal distribution is one with a clear peak and includes distributions, such as Gaussian, uniform, and chi-square distributions. For details, please see~\cite{stellato:Thesis2014}.}. 
An illustration of the uncertain estimation of a VESS's energy capacity and SOC is illustrated in Fig.~\ref{fig:Level1MPCproblems}. Note that in this manuscript, we make no assumption on the specific distribution of the estimation errors. However, for simulation purposes, the estimation errors are realized as normally distributed random variables to simplify simulation setup. 

\begin{figure} [h!]
    \centering
  \includegraphics[width=0.9\columnwidth]{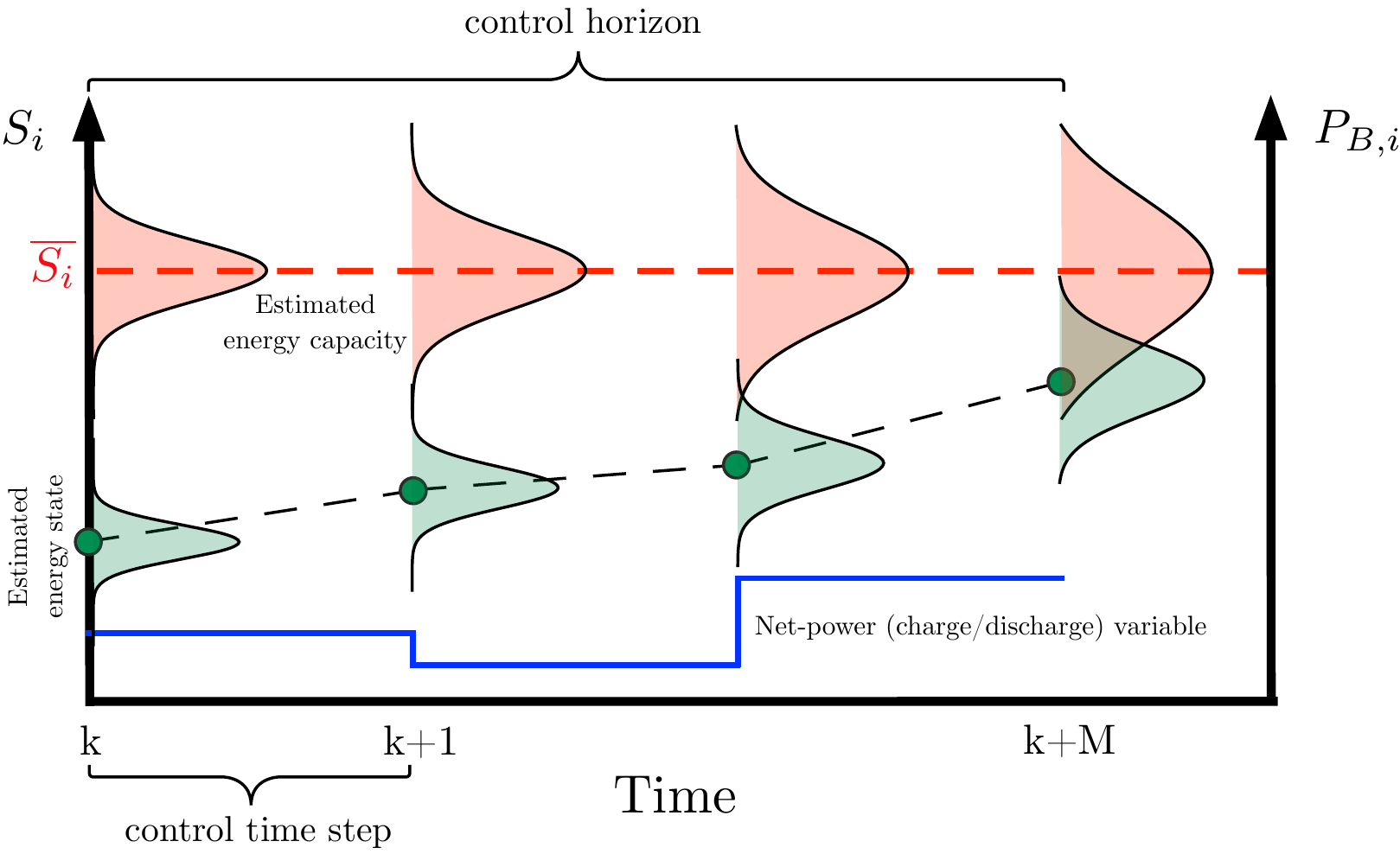}
  \caption{{\small Uncertainty in estimation of VESSs energy capacity (red) and initial SOC (green). The variance of uncertainties increases over time as the distance from current time step increases. The mean of the red distribution may also change randomily thru the prediction horizon, as long as it is predicted well.}} 
  \label{fig:Level1MPCproblems}
\end{figure}

\begin{definition}[Dynamic capacity saturation (DCS)]
Charging and discharging commands of the uncertain VESSs can be optimized with just the first moment, which we  denote as the Deterministic case. However, the underlying, uncertain VESS energy capacity may realize itself unexpectedly and saturate in the energy state, which zeros out the charging rate of the optimized control (power) action. We call this saturation phenomenon dynamic capacity saturation (DCS)~\cite{amini2018tradingoff, Hao:VB2015}. DCS may lead to unexpected power imbalances in the power system. To regulate these unexpected DCS-induced power imbalances, grid operators must rely on (expensive) generation to supply the difference based on their participation factor $d_i$, as shown below:
\begin{align}
\Delta P_{Gi}[k+1] &= - d_i\sum_{i \in N_{b}} \max \left\{\frac{S_i[l+1]-\overline{S_i}}{T_s}, 0\right \}. \label{eq:compens}
\end{align}
\end{definition}

\subsection{Chance constrained formulation}
Chance constrained optimization is employed to reduce the risk of DCS and to solve an optimization problem with uncertain parameters. The chance constraints should be satisfied with a  predefined probability level $1-\epsilon$, where $\epsilon \in (0,1)$ is the acceptable-worst-case violation level. Reducing risk increases system reliability and operational cost, which implies a clear trade-off. Within the context of pay-for-performance ancillary services~\cite{isone}, the operational costs are defined herein by the generators' reference-tracking errors. 

Chance constrained optimization problems can be solved with a probabilistically robust scheme, inspired by the so-called scenario approach. In the scenario approach, the chance constraint is substituted with a finite, but large number of deterministic constraints corresponding to different realizations of the underlying uncertainty space~\cite{vrakopoulou2013probabilistic}. By employing an adequate number of scenarios from this set (i.e. $N >> 1$), the approach is able to provide \textit{a-priori} guarantees of satisfying the chance constraint.
The scenario-based approach is useful in offline planning studies when the uncertainty is complex and captured via historical data, as it makes no assumption on the underlying distribution of the uncertainty. However, the number of scenarios required is a function of $\epsilon$ and the number of uncertain parameters can grow very large. If the underlying problem is convex,
there exists techniques to reduce the number of scenarios and mitigate computation by reformulating the problem into a robust optimization problem~\cite{Margellos2014}.

Indeed, if an accurate analytical model of the uncertainty distribution is known, the method \textit{analytical reformulation} can be employed to transform the chance constraint into a robust, deterministic constraint~\cite{roald2017corrective}. In contrast to the scenario approach, the analytical reformulation does not require sampling complex distributions or large data-set. This means that only a single reformulation for each chance constraint is needed, which makes implementation tractable at the timescale of interest~\cite{li2015analytical}. 

Next, we introduce the chance constraints related to the uncertain variables of the VESS (i.e., energy capacity and SOC) and briefly describe the analytical reformulation to derive a convex program. The formulation is presented with respect to the upper energy capacity limit of the VESS, but the lower limit can be handled in a similar manner. Note that unlike much of the literature of analytical reformulation in power systems, the work herein focuses on uncertain, controllable (time-coupled) energy resources rather than power injections.

\subsection{Analytical reformulation of chance constrained problem} 

Recall, the evolution of the SOC of the $i^{\text{th}}$ VESS over the prediction horizon (i.e., $l \in [1, M]$), is related to the estimated SOC of the $i^{th}$ VESS at time~$k$ (i.e., $S_{k,i}^{\text{est}}$) and charging (discharging) control actions as follows
\begin{align}
  &  S_i[l] = S_{k,i} + \sum_{m = 1  }^{l} \Delta S_i[m] 
\end{align}
where $\Delta S_i[l] := \Delta T_s(\eta_{\text{ch},i} P_{\text{ch},i}[l] - \eta_{\text{dis},i}^{-1} P_{\text{dis},i}[l])$ is the change of SOC at each timestep due to charging or discharging actions.

At time $k$, the estimated SOC of VESS $i$ is assumed to be a random variable centered on its true mean (i.e., $S_{k}^{\text{est}} = S_{k,i}^{\text{act}}+ \xi_{s,i}$), where the SOC estimation errors denoted by $\xi_s \in \mathds{R}^{N_B}$, with $\mu_s \in \mathds{R}^{N_B}$ as its mean and $\delta_s \in \mathds{R}^{N_B}$ as its standard deviation.  Any VESS technology that does not directly measure each DER's energy state frequently requires a dynamic state estimator that is specific to the model and information exchanges that underpin each specific VESS. However, this work makes no assumption on specific estimation and prediction methods. It is only assumed that these methods are uncertain and that the mean and variance of the initial SOC estimation errors can be estimated by offline or online data-driven approaches. In addition, due to the i.i.d. nature of DERs' end usage, estimation of the energy capacity of the $i^{\text{th}}$ VESSs can be modeled as a random variable around its true mean (i.e., $\overline{S}_{\text{est},i} = \overline{S}_{\text{act},i}+ \xi_{c,i}$), where $\xi_c \in \mathds{R}^{N_B}$ denotes the VESSs' capacity estimation error, with $\mu_c \in \mathds{R}^{N_B}$ as its mean and $\delta_c \in \mathds{R}^{N_B}$ the standard deviation, which can be determined from data-driven information scheme as in~\cite{zhang:aggDRuncertain2016}.

Then, for $\forall i = 1,\hdots,N_\text{B}$ the following constraints are equivalent:
\begin{align}
& \mathds{P}\left ( S_{k,i}^{\text{act}} + \sum_{m = 1}^{l} \Delta S_i[m] - \overline{S}_{\text{act},i}  \leq 0 \right) \geq 1-\epsilon \label{eq:prob0} \\ 
  &  \mathds{P}\left( S_{k,i}^{\text{est}} - \xi_{s,i} + \sum_{m = 1  }^{l} \Delta S_i[m] - \overline{S}_{\text{est},i} + \xi_{c,i} \leq 0 \right) \geq 1-\epsilon \label{eq:prob} \\
 &   S_{k,i}^{\text{est}} + \sum_{m = 1  }^{l} \Delta S_i[m] \leq \overline{S}_{\text{est},i} - f_{i,\epsilon}^{-1} \sqrt{\delta_{s,i}^2 + \delta_{c,i}^2 + \rho^{cs}_i\delta_{c,i}\delta_{s,i}}, &  \label{eq:robust}
\end{align} 
{\color{black}
where $f^{-1}_{i,\epsilon}$ denotes a \textit{safety-factor function} and $\rho^{cs}_i$ is the correlation coefficient of capacity and initial SOC estimation error for VESS $i$~\cite{SUMMERS2015116}. Thus, robustness against the uncertainties naturally begets an \textit{uncertainty margin} that is the product of the safety-factor function and the the variances and defines how much the constraint is tightened. The safety-factor function is defined by the information available on the underlying distribution of the uncertainty. For example, if we know that the uncertainty is normally distributed, then the safety-factor function is given by the (less conservative) inverse cumulative density function (cdf) of the standard Gaussian distribution. If we do not know the exact distribution, but have evidence that it is unimodal, the safety-factor function can be computed numerically with the Chebyshev generating function (CGF)~\cite{stellato:Thesis2014}, but is more conservative.  Finally, if we can only estimate first and second moments of the distribution and have no other information about the VESS's uncertainty, then we have to be robust against arbitrary distributions, which can also be achieved with the CGF~\cite{SUMMERS2015116}. The CGF allows us to take into account the severity of potential violations, which produces a very conservative outcome. This is expected since we do not assume any information about the distribution and still want to be robust against $\epsilon$-worst cases. 

Figure~\ref{fig:safetyFactor} illustrates the safety-factor functions for these different information scenarios. Clearly, knowing the exact distribution or even information that prescribes a unimodal distribuution leads to a reduction in conservativeness compared with only knowing the first two moments.
\begin{figure}[h]
      \centering
   \includegraphics[width=1\columnwidth]{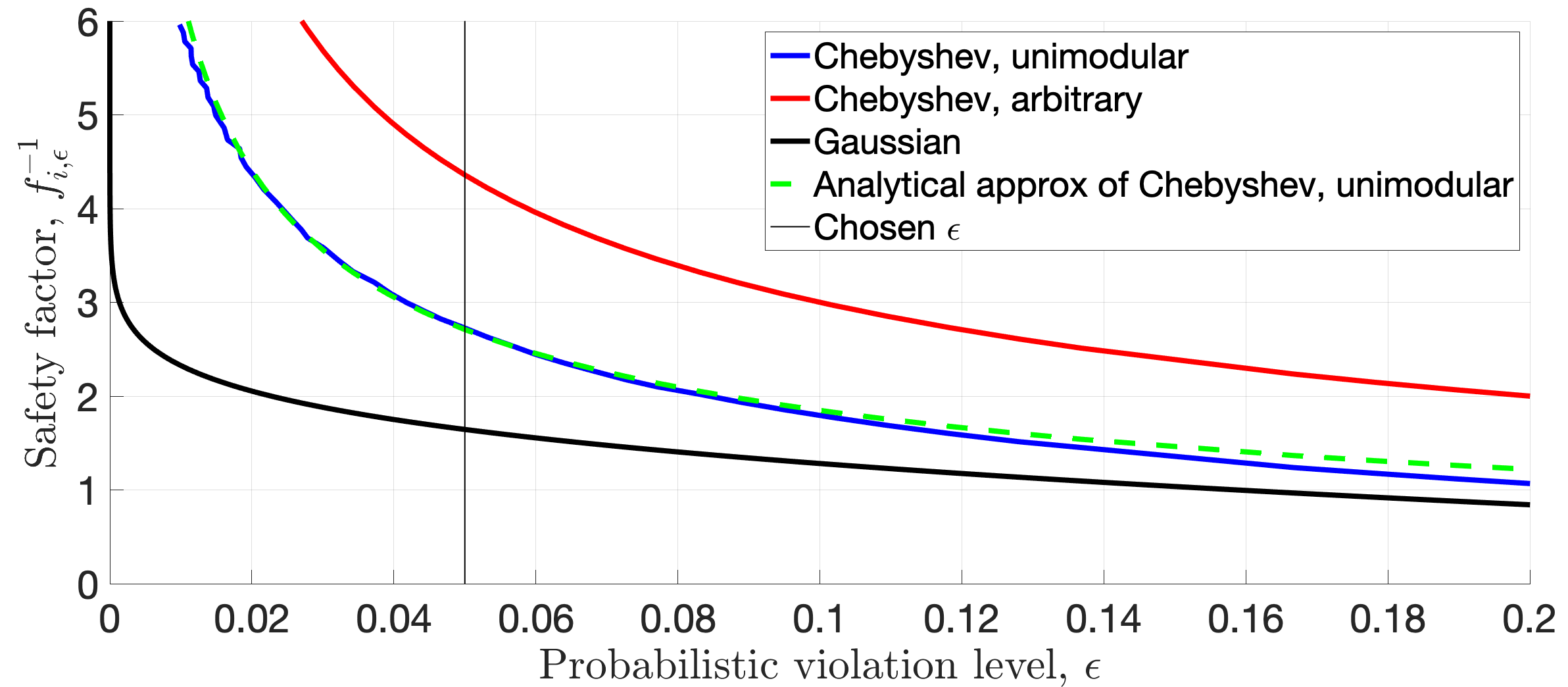}
      \caption{Comparing the effects of having different information available about the VESS on the VESS's safety-factor function. Larger values imply more conservative use of available flexibility. The results in this paper are based on choosing $\epsilon=0.05$.}
      \label{fig:safetyFactor}
   \end{figure}
The specific $f^{-1}_{i,\epsilon}$ for each of these three scenarios can be analytically or numerically\footnote{The safety-factor function for the unimodular distribution presented herein is a simple analytical approximation based on the exact numerical solution  from~\cite{stellato:Thesis2014}.} computed as:
\begin{align}
 f^{-1}_{i,\epsilon}  := 
 \begin{cases}
    \sqrt{\frac{1-\epsilon}{\epsilon}},& \text{for arbitrary distribution}\\
    \left(\frac{1-\epsilon}{e\epsilon} \right)^{1/1.95}, & \text{for unimodular distribution} \\
    \sqrt{2} \text{erf}^{-1}(1-2\epsilon), & \text{for Gaussian distribution}
\end{cases}
\end{align}
 The unimodular approximation shown in Fig.~\ref{fig:safetyFactor} is based on fitting a simple function (dashed green line) to the exact numerical values from~\cite{stellato:Thesis2014} (solid blue line). The approximation results in a simplified, closed-form safety-factor function with a relative error less than 15\% for $\epsilon < 0.20$ (and less than 5\% for $\epsilon < 0.10$) and is an inner approximation (i.e., more conservative). 

}

\begin{remark}
With an aggregation of flexible loads, knowing the exact parameters of the VESS (e.g., state of charge and energy capacity) by measuring each DER separately is challenging in real-time once coordinating DERs at scale. Thus, we do not assume we know these values; only that we can estimate the mean and variance of the parameters.  These estimates could come from data analytics that are re-run every couple of hours or days to update/true-up predictive models of the energy capacity and SOC estimators. Indeed, sparse live measurements available to the VESS may improve the estimates of the underlying distribution online (e.g., detect unimodal conditions and reduce moment uncertainty~\cite{stellato:Thesis2014}). That is, sparse live data streams may play an active role in updating the VESS parameters  to reduce conservativeness of the robust and risk-based schemes~\cite{Guo2018P2,zhang:aggDRuncertain2016}.
\end{remark}


\section{Risk-based Approach}


For a given chance constraint $\mathds{P}(f(x,\xi) \leq b) \geq 1 - \epsilon$, where $x$ and $\xi$ are decision and uncertain variables, the magnitude of constraint violation is a function of $\xi$ and is given by $y(\xi) = f(x,\xi) - b $. Negative $y$ indicates constraint satisfaction while positive $y$ implies constraint violation. The chance constraint limits the probability of violation ($y>0$) to a predefined risk limit, $\epsilon$. Since chance-constrained approaches ignore the severity of the constraint violation, the approaches are conservative and a closed-loop chance constrained MPC (CC-MPC) implementation may lead to reduced performance of the system (by significantly reducing the   available flexibility). Authors in~\cite{roald2017optimized}, consider the severity of a constraint violation, by weighting the probability of the constraint violation by the magnitude of the constraint violation.

In applying chance constraints, there is a clear trade-off between high reliability (i.e., conservative uncertainty margin for VESSs) and low nominal cost (i.e., use as much VESS as possible), which depends on how risk limits are chosen. Risk limits are generally chosen as a predefined parameters (e.g. $\epsilon \in (0.90 \quad 0.99)$) based on the importance of the constraint.

Unlike the robust approach that limits the SOC of VESSs to a predefined robust limit $\overline{S}_{\text{rob}}$, we propose a novel risk-based approach that allows the solution to exceed the robust limit at each point in time. This is possible by introducing the operating risk, $\mathcal{R}$, which is a new decision variable. Thus, performance and risks can be co-optimized, which leads to the following multi-objective optimization problem:

\small

\begin{subequations}
  \begin{align}
    &  \underset{\textstyle P_{G},P_{{\text{ch}}},P_{{\text{dis}}},\mathcal{R}}{\text{min}} \quad 
      \sum_{l = k}^{k+M} J_1[l] +J_2[l] \hspace{90pt} \label{eq:sec2opt} \\
   & \text{s.t.} \quad  \text{\eqref{eq:const2} to \eqref{eq:const19}, }  \nonumber \\
   &  \mathcal{R}_i[l] = \text{max}(0,S_i[l]-\overline{S}_{\text{rob},i})  \label{eq:max1},  \\
   &   \mathcal{R}_i[l]\le \overline{\mathcal{R}_i}, \\
   & \overline{\mathcal{R}_i} = \overline{S}_{\text{est},i} - \overline{S}_{\text{rob},i}
  \end{align}
\end{subequations}
\begin{align}
    & \text{where}  \nonumber  \\
     & J_1[l] :=  \sum_{\forall i \in N_g} c_{G,i} ( P_{G,i}[l]-P_{G,i}^\text{ref}[l] )^2+\sum_{\forall ij \in \mathcal{E}}c_{T,ij} {(\Delta\hat{T}_{ij}[l]})^2 \nonumber \\
   & J_2[l] :=  \sum_{\forall i \in N_B} c_{\mathcal{R},i}(\mathcal{R}_i[l])^2 \nonumber
\end{align}  
\normalsize

The robust limit can be computed by analytical reformulation or scenario based approach or determined by using expert knowledge. As the risk-based MPC receives updated estimates at each time step, the cost of risk, $c_{\mathcal{R}}$, can be designed such that it penalizes  risk early in the horizon and lowers the penalty later in the horizon. Larger $c_{\mathcal{R}}$ indicates higher cost and higher security and it is necessary to reach a good balance between risk of a VESS DCS and nominal tracking performance. To relate the value of improving tracking performance and the associated increase in operational risk, an efficient frontier for the tracking performance versus operational risk can be computed~\cite{deb2016multi}.

The risk-based MPC seeks to drive the SOC of VESSs below the robust limit, but once below the robust limit, there is no incentive to further lower the SOC, as shown in~\eqref{eq:max1}. This constraint reformulated as a set of linear constraints: $0 \leq \mathcal{R}_i[l]$ and $S_i[l]-\overline{S}_{\text{rob},i} \leq \mathcal{R}_i[l]$. Figure~\ref{fig:Guass2} illustrates an example of the SOC of the VESS with respect to the estimated capacity of the VESS and its robust limit. 

Thus, the three approaches developed herein are: deterministic, robust, and risk-based. They are summarized as follows:
\begin{enumerate}[label=\Roman*.]
    \item The deterministic method dispatches VESSs with respect to just the first moment (averages) for the estimated SOC and predicted energy capacity.
    \item The robust method dispatches VESSs with respect to robust limits calculated with analytic reformulation (illustrated with the unimodular and Gaussian distributional assumptions).
    \item The risk-based chance constrained (RB-CC) method co-optimizes reference-tracking performance and operational risk of DCS. Note that by sweeping $c_{\mathcal{R}}$ from $0 \rightarrow \infty$, the performance of the controller changes from the deterministic approach to the robust approach.
\end{enumerate}

\begin{figure}[t!]
      \centering
   \includegraphics[width=0.88\columnwidth]{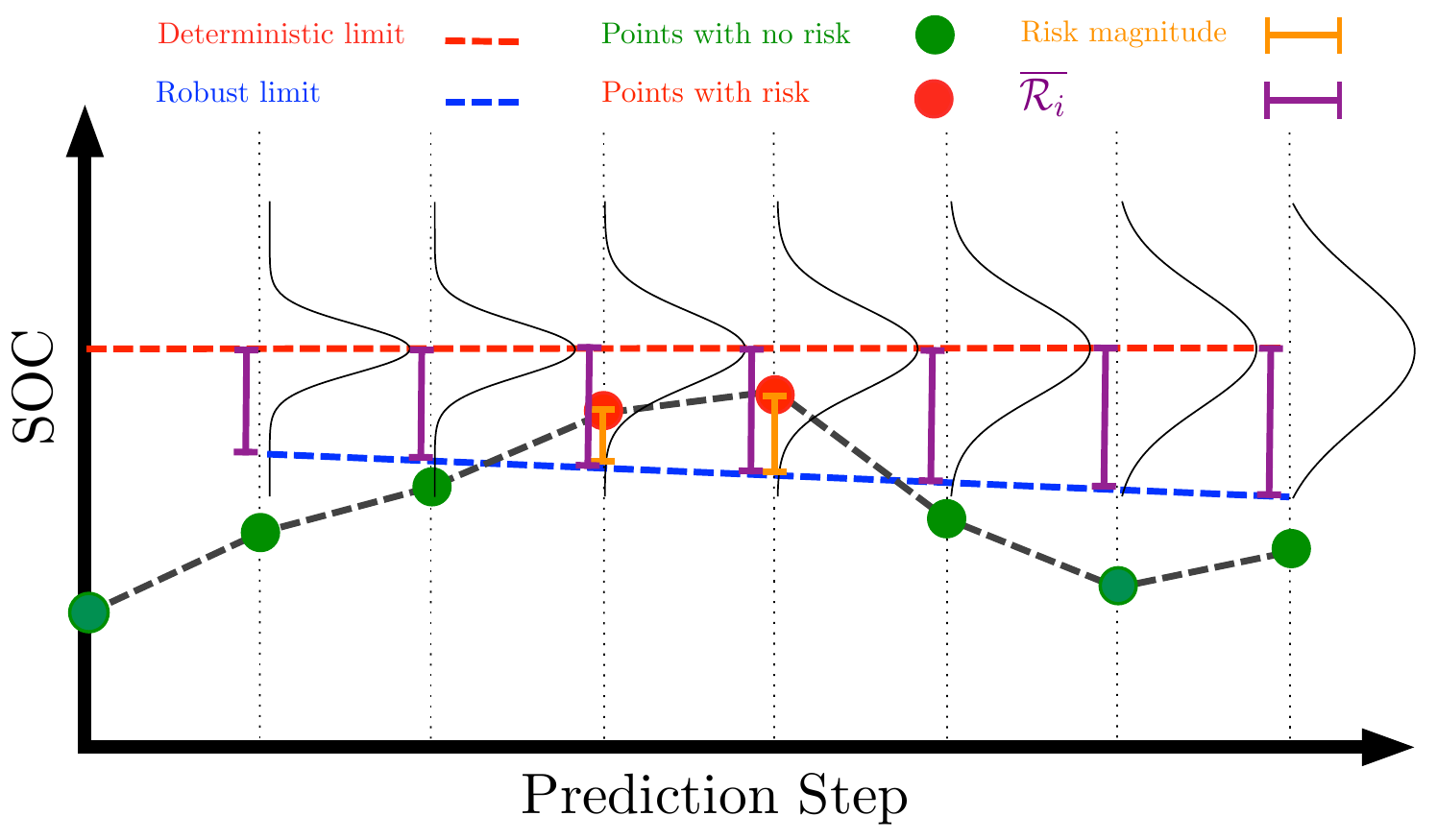}
      \caption{ {\small Illustrating the evolution of the SOC of a VESS with respect to the estimated capacity and robust bound and the corresponding risk imposed on the system performance.The variance of uncertainties grow over time as the distance from current time step increases. The green circles highlight the points with zero added risk. On the other hand, the red circles demonstrate when the VESSs' SOC is greater than the robust limit and takes on an increased, but weighted risk of DCS.}}
      \label{fig:Guass2}
   \end{figure}
   
Recall that unlike the existing literature on chance constrained optimization in power systems, this work considers the uncertainty of  \textit{time-coupled energy variables} on a fast time scale for corrective control.

\section{Simulation results and discussions}
In this section, the proposed control scheme in Fig.~\ref{fig:overview} is demonstrated on an augmented version of the IEEE RTS-96 power system test case. The system is detailed in~\cite{grigg1999ieee}. 
All optimization problems are solved in MATLAB and AMPL using the solver GUROBI. MPC employs a simplified linear model to compute all optimal control actions. However, the  plant model is the AC network, with line temperature computed based on the non-linear thermodynamic IEEE Standard~738 conductor temperature model to accurately capture the effects of implementing the MPC actions.
The aim of this case study is to demonstrate generator reference-tracking performance with uncertain VESSs while considering physical constraints of the power system. Since no data regarding the actual estimation errors is available, the estimation errors are generated from Gaussian distribution, which satisfy the unimodular condition. Under this condition, we can compare the performance of the controller when the robust bound is chosen based on the unimodular (Chebyshev) assumption with no other knowledge of the distribution and the less  conservative Gaussian assumption. Since the IEEE RTS-96 system is designed as a highly reliable system with high thermal ratings for lines, nominal thermal ratings are reduced by $40\%$ to reduce line temperature limits to a range of $60-70 \celsius$. The network parameters are shown in Table~\ref{table:param}.


\begin{table}[]
    \centering
    \caption{Three large VESSs Simulation Parameters}
\setlength{\tabcolsep}{16pt}
\begin{tabular}{ l l}
\toprule
     Description & Value       \\ \midrule
     Number of buses & $73$ \\
     Number of branches & $120$ \\
     Number of generators & $96$ \\
     Total load & $8550$ MW \\
     Number of VESSs & $3$ \\
VESS bus ID (location) & $11, 35, 59$ \\
VESS energy capacity (MWh) & $250, 300, 350$ \\
VESS initial SOC(\%) & $40,50,60\%$ \\
Maximum VESS power (MW) & $250, 300, 350$  \\
VESS ramp rate limits & $60$ MW/min \\
     Sampling time & $60$ sec \\
MPC prediction horizon  & $20$ mins \\
Violation level $\epsilon$ & 0.05 \\
Avg. MPC solve time & $3.19$ sec \\
\bottomrule
\end{tabular}
\label{table:param}
\end{table}

\begin{table}[]
    \centering
    \caption{VESS Uncertainty parameters}
\setlength{\tabcolsep}{20pt}
\begin{tabular}{lcc}
\toprule
  \multicolumn{1}{l}{Case}     & \multicolumn{1}{l}{$\delta_{s,i}$ (\%)} & \multicolumn{1}{l}{$\delta_{c,i}$ (\%)}\\
       \midrule
\multicolumn{1}{l}{\begin{tabular}[l]{@{}c@{}} Low uncertainty \end{tabular}} &  2 & 3 \\ 
\multicolumn{1}{l}{\begin{tabular}[l]{@{}c@{}}High uncertainty \end{tabular}} &  5 & 10 \\ 
\bottomrule
\end{tabular}
\label{tbl:CaseABCD}
\end{table}


Initially, the system is at steady-state (i.e., generators following exactly an economic trajectory and VESS resources being available for balancing reserves), but at time $t=5$ mins, the system experiences a net-load disturbance (e.g., forecasted net-load) that requires VESS balancing reserves to provide 10\% decrease in the net-load (i.e., $855$ MW) to minimize unnecessary generators ramping. In order to have a realistic case, VESSs are sized and initialized differently as shown in Table~\ref{table:param}} to compare the performance of the proposed RB-CC method aginst deterministic and robust approaches, 
$N_T=100$~trials (i.e., realizations) are performed for two different scenarios of uncertainty as shown in Table~\ref{tbl:CaseABCD}. \textcolor{black}{Since the baseline consumption and consequently the capacity of VESSs are dependent on the same types of uncertain parameters, capacity estimation errors of VESSs are assumed to be correlated.} The sum of squared error (SSE) of reference tracking of generators over the entire simulation time
\begin{align}
   \text{J}_{\text{Gen}} :=  \sum_{l=1}^N \sum_{i\in N_g} (P_{G,i}[l] -   P_{G,i}^\text{ref}[l])^2
\end{align}
is used as the tracking MPC performance metric.

\begin{figure}[ht!]
      \centering
   \includegraphics[width=0.88\columnwidth]{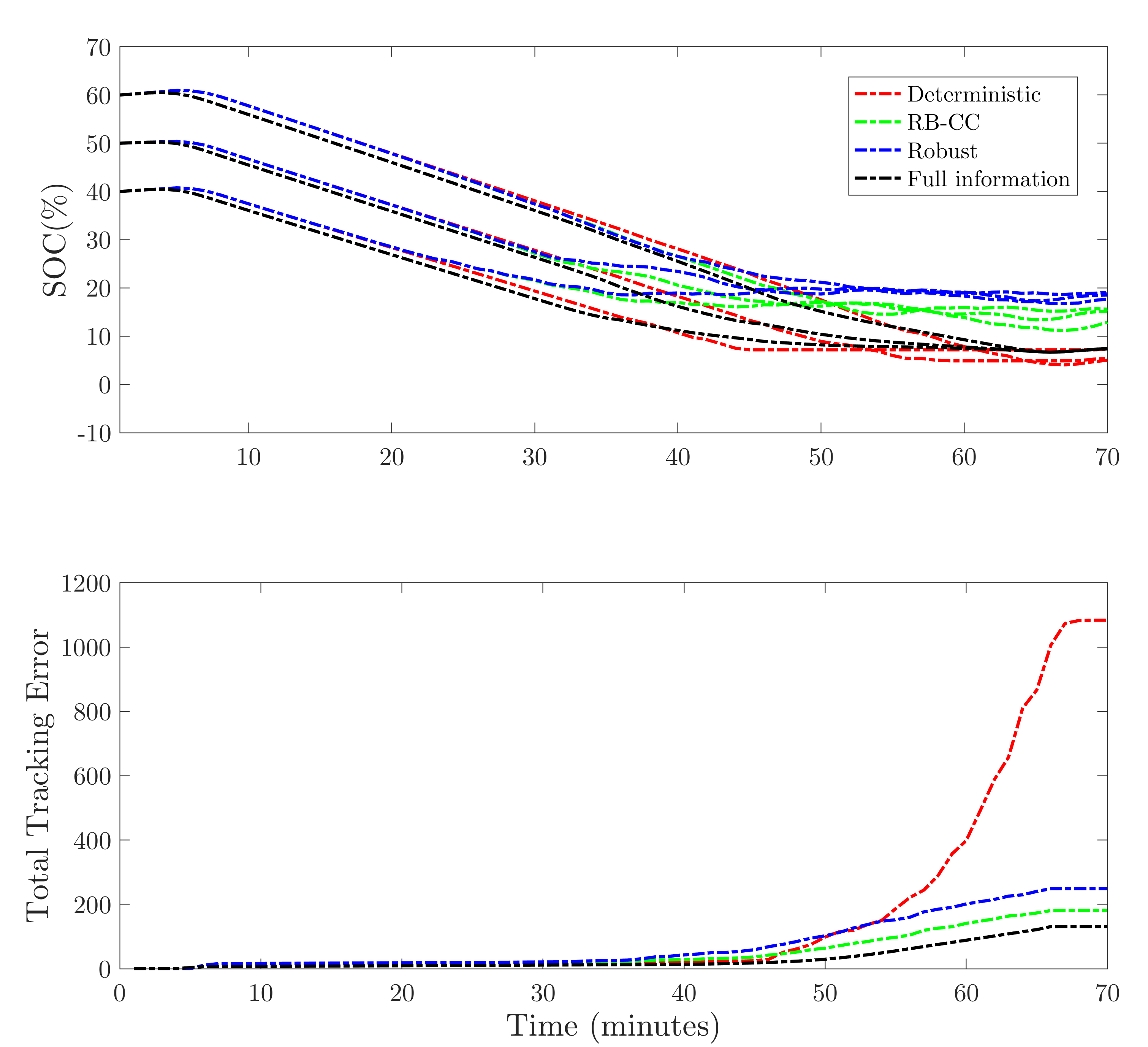}
      \caption{{\small Cumulative squared tracking error (bottom plot) and evolution of the energy state of the charge of the VESSs (top plot) under full information, deterministic, robust, RB-CC approaches for one randomly chosen trial.}}
      \label{fig:timeseries}
   \end{figure}

To better understand the role of uncertainty on the performance of the system, the MPC problem is solved under the assumption that the true capacity and state of the charges of the VESSs are available to the controller (full or perfect information) and shown in Fig.~\ref{fig:timeseries}. This provides a baseline against which to compare the tracking performance of deterministic, robust, and RB-CC approaches. Note that even under the full information assumption, since the VESSs are energy-constrained resources, the tracking error is not zero. That is, the optimal solution involves regulating both generators and VESS resources against ramp-rate limits.
Mean and standard deviation of the reference tracking of generators $\text{J}_{\text{Gen}}$, for full-information, deterministic, robust, and RB-CC approaches under the different scenarios of uncertainty for (unimodal) Chebyshev and Gaussian approaches are shown in Tables~\ref{table:compare}. The robust approach is used as a benchmark to evaluate the proficiency of the proposed RB-CC approach. Smaller SSE, implies lower tracking error. Poor performance of the deterministic approach is due to dispatching VESSs without considering the second moment of the uncertainty which increase the risk of DCS. By employing the robust approach, chances of DCS are low, but the method is overly-conservative, which under-utilizes the VESS flexibility. However, with the proposed RB-CC approach, the controller is able to optimally trade off the flexibility available from the VESSs while explicitly considering the uncertainty by accepted increased risk when it is most valuable to do so. Therefore, with both (unimodal) Chebyshev and Gaussian methods, the RB-CC outperforms the corresponding robust approaches. Note that the Gaussian method outperforms the unimodal Chebyshev method (as expected) since the assumption of a Gaussian distribution is much stronger than a unimodular assumption. Also, since the unimodal Chebyshev assumption yields a very conservative robust bound, the RB-CC approach can achieve significant improvements by taking (a calculated) risk. In fact, for low uncertainty, the RB-CC (unimodal Chebyshev) has similar performance to Robust (Gaussian).

\renewcommand{\tabcolsep}{6pt}
\begin{table}[ht!]
\caption{Comparing mean and standard deviation of tracking performance of the MPC (i.e., $\text{J}_{\text{Gen}}$) under full-information, deterministic, robust, and RB-CC approaches}
\label{table:compare}
\begin{subtable}[ht]{\linewidth}
\resizebox{\columnwidth}{!}{%
\begin{tabular}{lcccccc}
\toprule
 & \multicolumn{3}{c}{Low uncertainty} & \multicolumn{3}{c}{High uncertainty} \\ [0.2 cm] \hline
 & $\mu$ & $\sigma$ & $\mu / \mu_{\text{Rob}} (\%)$ & $\mu$ & $\sigma$ & $\mu / \mu_{\text{Rob}} (\%)$ \\ [0.1 cm]\hline 
Full-information & 82.4 & 31.9 & 43.6 & 102.4 & 56.5 & 29.7\\
Deterministic & 367.1 & 622.5 & 194.0 & 1239.3 & 1796.1 & 359.4 \\
Robust & 189.2 & 33.6 & 100.0 & 344.8 & 182.4 & 100.0 \\
RB-CC & 166.6 & 194.6 & 88.1 & 304.1 & 457.2 & 88.2\\
\bottomrule
\end{tabular}%
}
\caption{Robust bound is determined with specific knowledge of estimation error distribution (Gaussian approach)}
  \end{subtable}
  \begin{subtable}[ht]{\linewidth}
\resizebox{\columnwidth}{!}{%
\begin{tabular}{lcccccc}
\toprule
 & \multicolumn{3}{c}{Low uncertainty} & \multicolumn{3}{c}{High uncertainty} \\ [0.2 cm] \hline
 & $\mu$ & $\sigma$ & $\mu / \mu_{\text{Rob}} (\%)$ & $\mu$ & $\sigma$ & $\mu / \mu_{\text{Rob}} (\%)$ \\ [0.1 cm]\hline 
Full-information & 82.4 & 31.9 & 30.8 & 102.4 & 56.5 & 20.2\\
Deterministic & 367.1 & 622.5 & 137.2 & 1239.3 & 1796.1 & 244.6 \\
Robust & 267.6 & 18.6 & 100.0 & 506.7 & 36.4 & 100.0 \\
RB-CC & 188.9 & 10.5 & 70.6 & 368.7 & 42.8 & 72.8\\
\bottomrule
\end{tabular}%
}
\caption{Robust bound is determined without specific knowledge of estimation error distribution (unimodular Chebyshev approach)}
 \end{subtable}
\end{table}


   
\begin{figure}[ht!]
            \centering
      \begin{subfigure}{.5\textwidth}
   \includegraphics[width=0.88\columnwidth]{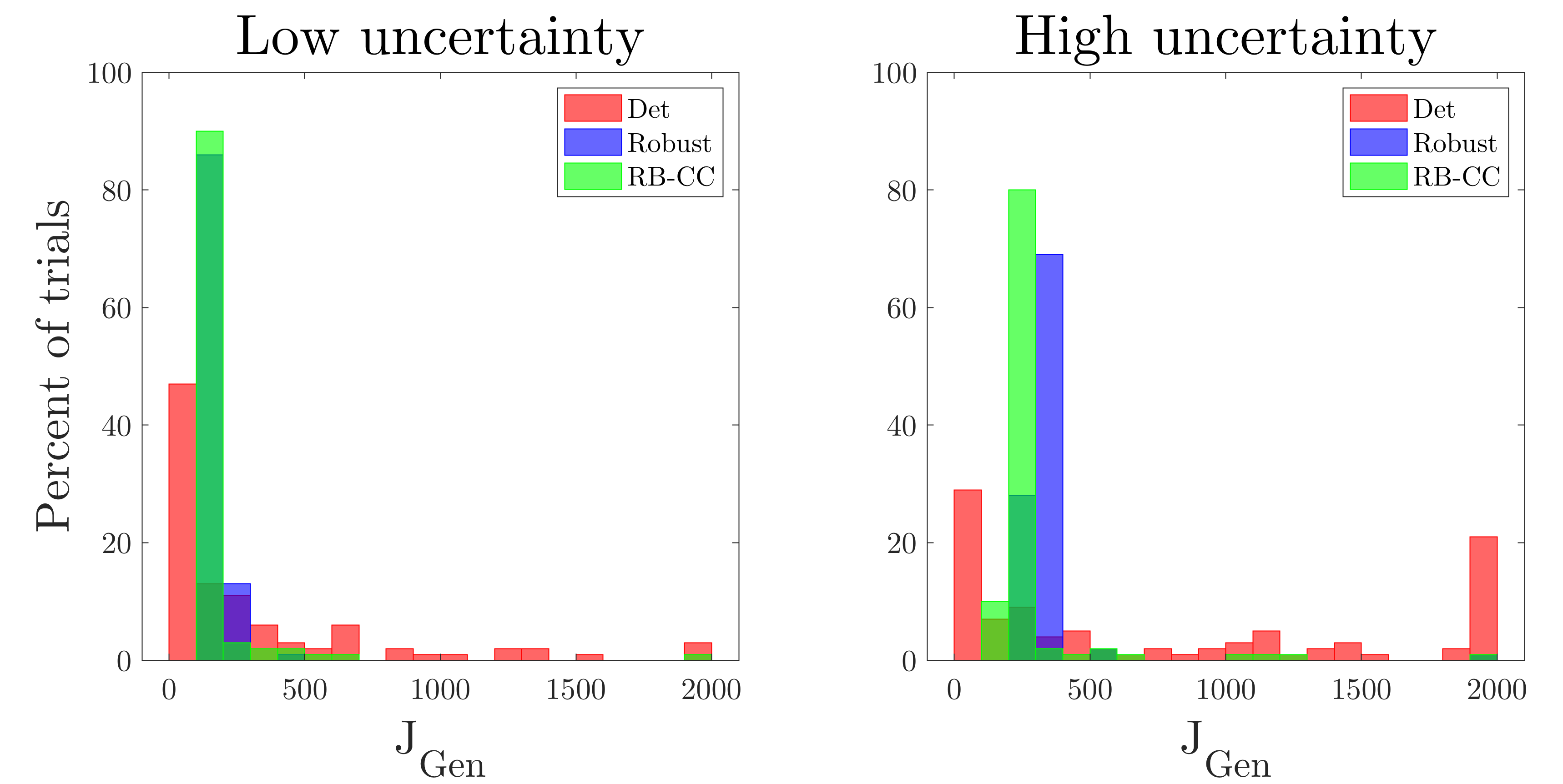}
     \caption{{\small Robust bound is determined with knowledge of estimation error distribution (Gaussian approach)}}
   \end{subfigure}
    \begin{subfigure}{.5\textwidth}
   \includegraphics[width=0.88\columnwidth]{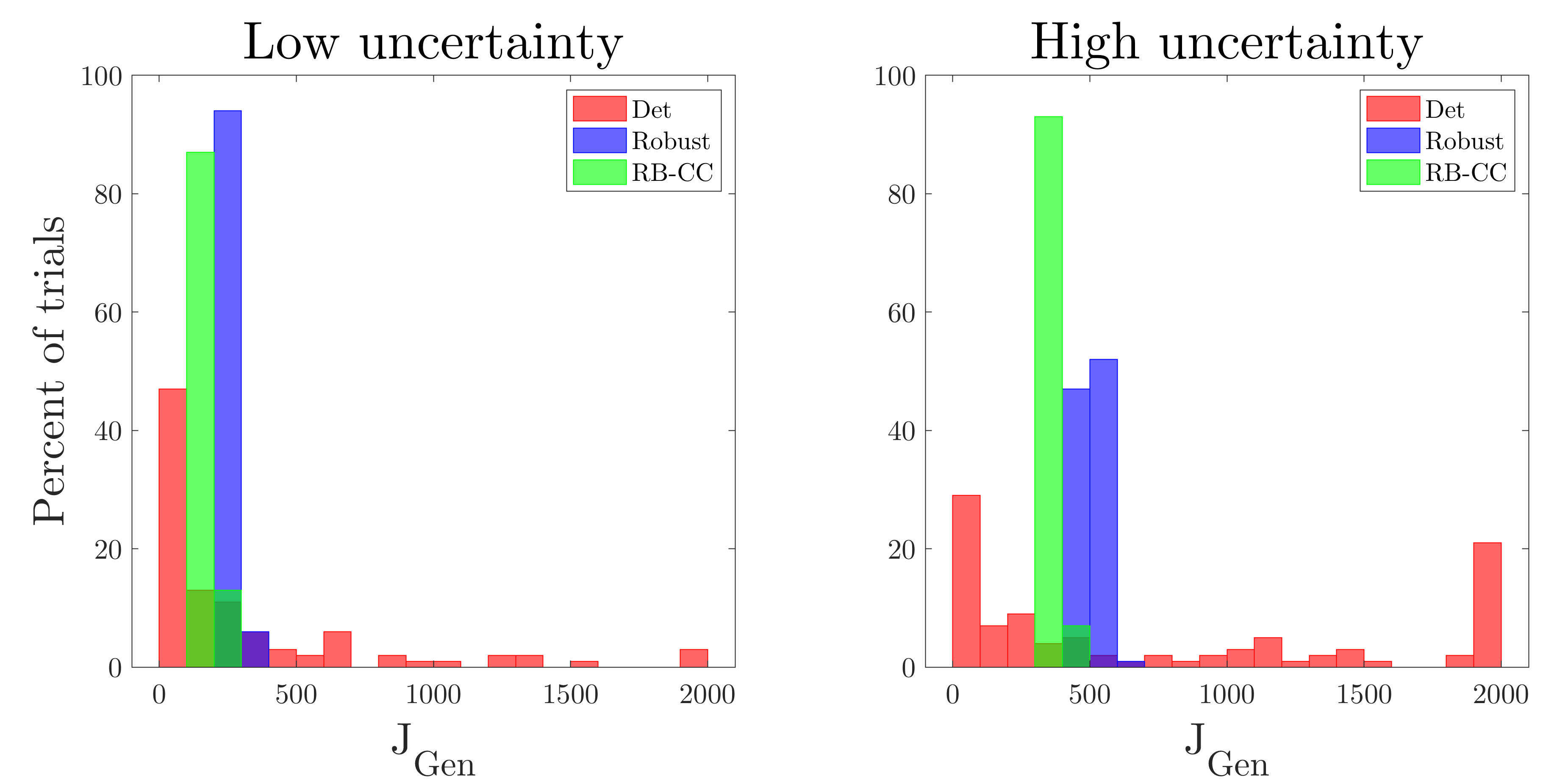}
     \caption{{\small Robust bound is determined with no knowledge of estimation error distribution (Chebyshev approach)}}
   \end{subfigure}
      \caption{{\small Histogram of the squared tracking error under deterministic, robust and RB-CC approaches. For visualization purposes, trials with squared tracking error of greater than $2000$ are categorized in the last bin.}}
      \label{fig:Hist4Case}
   \end{figure}

 Histograms for the total squared tracking error ($\text{J}_{\text{Gen}}$) for the deterministic, robust, and RB-CC cases for Gaussian and Chebyshev robust bound are shown in Fig.~\ref{fig:Hist4Case}. Note that, for both Gaussian and Chebyshev methods, under high levels of uncertainty, the deterministic formulation becomes susceptible to DCS, which reduces average closed-loop tracking performance while the RB-CC formulation outperforms the robust approach.

It is shown that regardless of how the robust bound is chosen (e.g., Chebyshev, Gaussian), the corresponding RB-CC outperforms the Robust approach. Therefore, for the rest of the paper, we focus on the unimodal Chebyshev method, which does not need exact knowledge of the underlying VESS distributions.

{\color{black}Figure~\ref{fig:timeseries} shows the cumulative generator squared tracking error and time-series evolution of the SOC of the VESSs for one of the 100 trials. The robust approach dispatches VESSs considering their robust SOC limit, which results in reduced tracking error performance at the beginning compare to the deterministic approach. However, as time goes on and the SOC of VESSs approach their limits, the chance of DCS and, consequently, the need for corrective actions increases, which results in reduced tracking performance for the deterministic case. Since the RB-CC approach explicitly considers the risk-performance trade off in the optimization, DCS is prevented and  high tracking performance is achieved.}

To further investigate the effectiveness of the proposed method, the three large VESSs (one located in each region), are replaced by nine smaller VESSs (three located in each region) and location and parameters of small VESSs are shown in Table~\ref{table:paramSVESS}. To have a realistic case, VESSs are designed to have different capacity and initial SOC. Intuitively, smaller VESSs should reduce the severity of DCS events, but increase their frequency. The same analysis has been carried out on the system with the small VESSs under high uncertainty scenario and as results are shown in Fig.~$7$ and Table~V, the RB-CC outperforms the robust approach significantly. 

\begin{figure}[ht!]
      \centering
   \includegraphics[width=0.88\columnwidth]{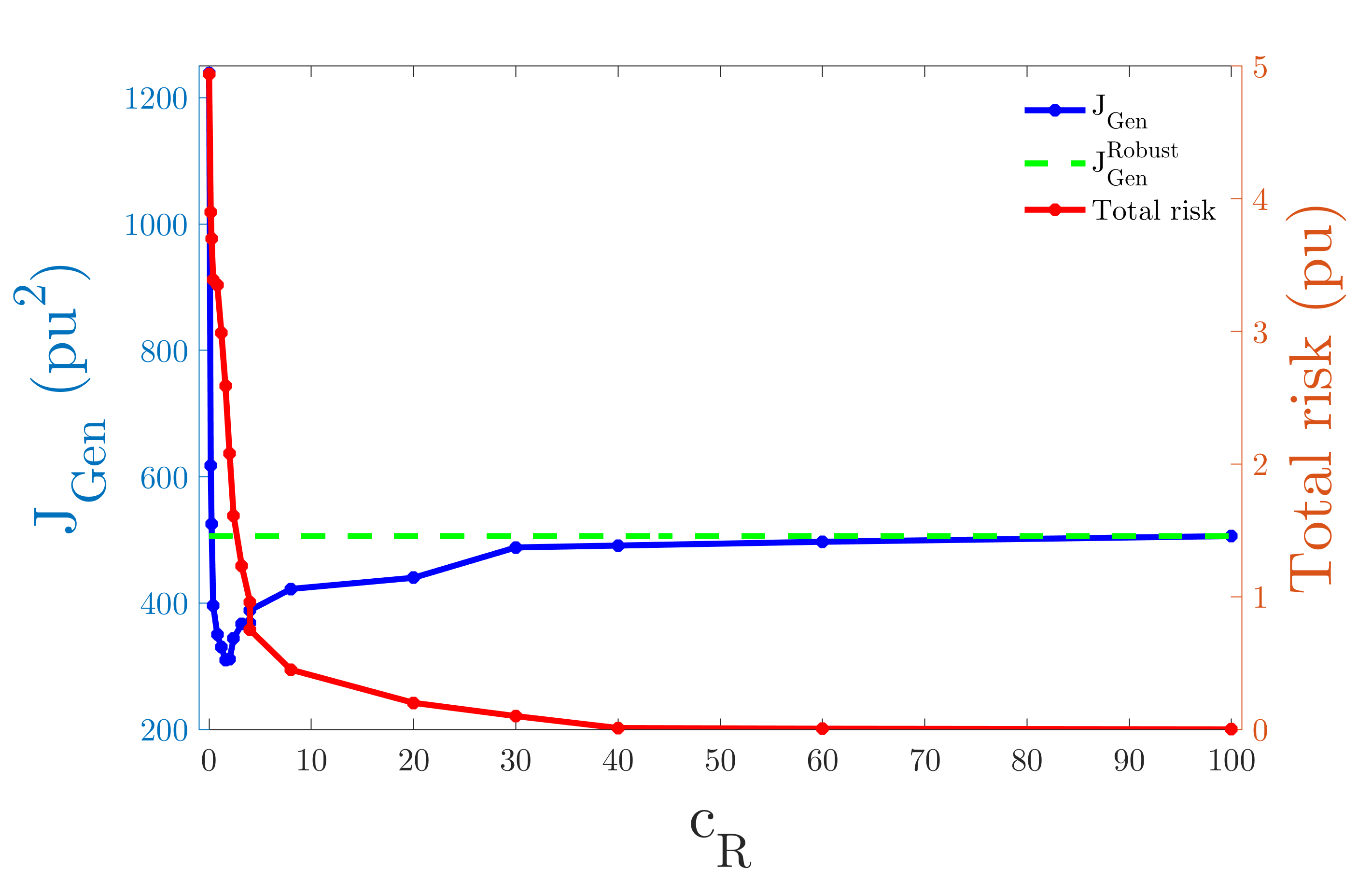}
      \caption{ {\small Role of the risk cost, $c_\mathcal{R}$ on the tracking performance and the total risk imposed to the power system.}}
      \label{fig:frontier}
   \end{figure}

{\color{black}To illustrate the role of risk cost, $c_\mathcal{R}$, on the tracking performance and total risk imposed to the power system, we scale the risk cost, $c_{\mathcal{R}}\in [0,100]$ where $c_{\mathcal{R}}\approx 0$ begets the Deterministic approach and $c_{\mathcal{R}}>40$ approximates the Robust approach.}  Figure~\ref{fig:frontier} shows the optimal trade-off between the tracking performance $J_{\text{Gen}}$, total operational risk that MPC accepts and risk cost. As expected, for small $c_{\mathcal{R}}$, the controller dispatches VESSs with respect to the expected energy capacity limit, disregarding the robust limit. This results in larger total risk, more DCS and poor tracking performance. For large $c_{\mathcal{R}}$, VESSs are dispatched conservatively and while the VESSs are dispatched at a very low risk level, parts of the flexibility offered by VESSs are declined.  

\begin{table}[ht!]
\centering
\caption{Parameters of nine small VESSs}
\label{table:paramSVESS}
\begin{tabular}{ l l }
\toprule
Description                     & Value       \\ 
\midrule
Number of VESSs                  & $9$ \\
Total energy capacity of VESSs   & $855$ MWh \\
Bus ID (location) of VESSs       & $11,17,24,35,41,48,59,65,72$ \\
VESS energy capacity (MWh)             & $95,85,100,85,90,95,80,100,90$ \\
Initial VESS state of charge(\%)     & $45, 65, 55, 50, 40, 60, 55, 45, 55$ \\
Maximum VESS power output (MW)        & $95,85,100,85,90,95,80,100,90$ \\
VESSs ramp rate limit            & $20$ MW/min \\
\bottomrule
\end{tabular}
\end{table}


\begin{table}
	\begin{minipage}{0.5\linewidth}
		\caption{Nine small VESSs under high uncertainty}
		\label{table:student}
		\centering
\resizebox{1\columnwidth}{!}{%
\begin{tabular}{lccc}
\toprule
 & $\mu$ & $\sigma$ & $\mu / \mu_{\text{Rob}} (\%)$ \\
 \midrule 
Full-information & 89.93 & 24.23 & 25.15 \\ 
Deterministic & 1327.3 & 965.0 & 371.2 \\
Robust & 357.6 & 20.1 & 100  \\
RB-CC & 204.3 & 110.0 & 57.1 \\
\bottomrule
\end{tabular}
}
	\end{minipage}\hfill
	\begin{minipage}{0.45\linewidth}
		\centering
		\includegraphics[width=40mm]{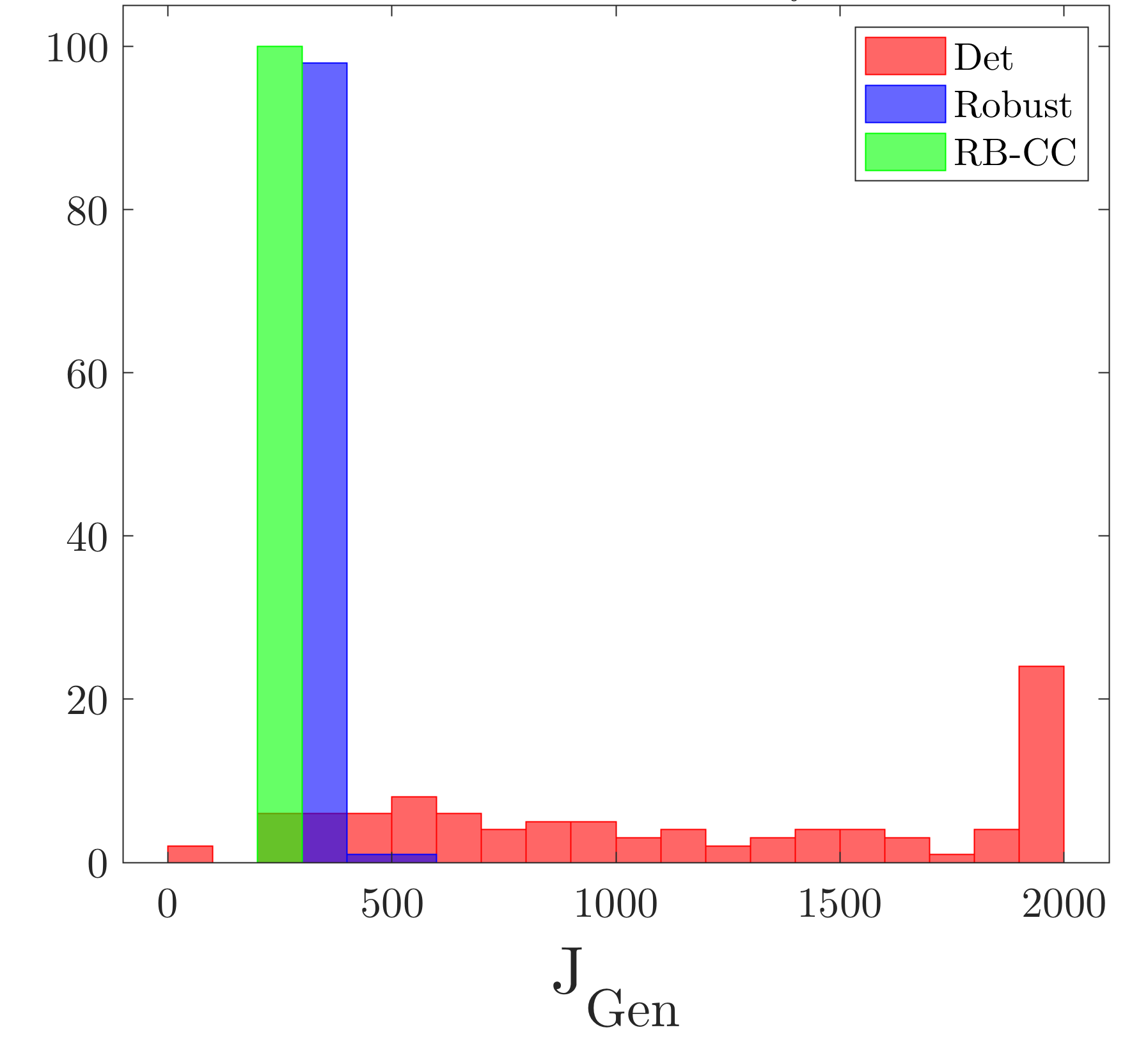}
		\label{fig:small_all}
	\end{minipage}
			\captionof{figure}{{\small Performances of the deterministic, robust and RB-CC method are analyzed in presence of nine VESSs (three in each region). (left) Comparing average tracking performance of deterministic, robust and RB-CC approaches. (right) Histogram of the squared tracking error under deterministic, robust and RB-CC approaches. For visualization purposes, trials with squared tracking error of greater than $2000$ are categorized in the last bin.  }}
\end{table}



\section{Conclusion and future work}
This paper studies the performance of a bi-level receding horizon predictive optimal power flow problem for managing variability with uncertain, flexible grid assets, such as VESSs. Since the SOC and capacity of VESSs can not be measured directly, a dynamic state estimator and simplified VESS aggregate model must be employed, which introduce uncertainty. This uncertainty in energy-constrained resources gives rise to the notion of dynamic capacity saturation (DCS). To overcome DCS, uncertainty can be managed by employing robust approaches. However, there is a sensitive trade-off between robustness of the optimized dispatch and closed-loop performance of the system. Indeed, robust approaches may lead to a conservative (high-cost) solution. Therefore, we introduced a RB-CC approach under which the operational risk is optimized with respect to the dynamic states of the VESSs over a receding horizon. The numerical studies indicate that RB-CC outperforms other methods and significantly reduces DCS while maintaining good tracking performance.


Future work will focus on further reducing conservativeness by augmenting the VESS uncertainty models with finite distributional support since we know that state-of-charge and energy capacities are non-negative entities~\cite{SUMMERS2015116}.  However, finite distributional support gives rise to additional linear matrix inequalities (LMIs), which turns the formulation into a semi-definite program (SDP), which is numerically more sensitive than what is presented herein. In addition, while the simulation results herein include a correlation between the VESSs' capacity estimation errors, we are interested in further generalizing sources of uncertainty and correlations between the VESSs. For example, the role of uncertain VESS power limit forecasts on corrective dispatch are of theoretical and practical interest. Secondly, we are interested in studying how managing uncertainty in the VESS with the risk-based approach affects the line-temperature regulation problem in the presence of disturbances from renewable generation. The problem formulation can already support this study, but the focus is different from what is presented so far. Finally, we are interested to explore linearized or convex AC formulations that can actively correct any possible voltage excursions to ensure that voltage profiles can be regulated effectively.

{\small
 \bibliographystyle{IEEEtran}
 \bibliography{ref}}

%



\end{document}